\documentclass[12pt]{amsart}

\usepackage[T1]{fontenc}
\usepackage[utf8]{inputenc}

\usepackage{times}

\usepackage{graphicx, xspace}
\usepackage[dvips]{hyperref}

\theoremstyle{plain}
\newtheorem{lemma}{Lemma}
\newtheorem{theorem}[lemma]{Theorem}
\newtheorem{corollary}[lemma]{Corollary}
\newtheorem{proposition}[lemma]{Proposition}

\theoremstyle{remark}

\newtheorem{example}[lemma]{Example}

\newcommand{\comment}[1]{}
\newcommand{\R}{\mathbb{R}}
\newcommand{\N}{\mathbb{N}}
\newcommand{\C}{\mathbb{C}}
\newcommand{\D}{\mathcal{D}}
\newcommand{\A}{\mathcal{A}}

\newcommand{\Grozciete}{\tilde{G}}
\newcommand{\Gobciete}{G'}
\newcommand{\GobcieteBardzo}{G''}

\newcommand{\E}{\mathbb{E}}
\newcommand{\Sy}{\tilde{S}}
\newcommand{\gwia}{^\ast}

\author{Valentin F\'eray}
\address{LaBRI, Universit\'e Bordeaux 1, 351 cours de la Lib\'eration, 33 400 Talence, France}
\email{feray@labri.fr}

\author{Piotr \'Sniady}
\address{Insitute of Mathematics, Polish Academy of Sciences, ul.~Śniadeckich~8,
P.O.B.\ 21, 00-956 Warszawa 10, Poland \newline 
 Institute of Mathematics,
 University of Wroclaw,  \mbox{pl.\ Grunwaldzki~2/4,} 50-384
 Wroclaw, Poland} 
\email{Piotr.Sniady@math.uni.wroc.pl}

\title[Asymptotics of characters of symmetric groups]%
{Asymptotics of characters of symmetric groups related to Stanley character formula}

\DeclareMathOperator{\Tr}{Tr}

\DeclareMathOperator{\vol}{vol}
\DeclareMathOperator{\supp}{supp}

\DeclareMathOperator{\inv}{cinv}
\DeclareMathOperator{\dimm}{dim}

\DeclareMathOperator{\orbits}{orbits}
\DeclareMathOperator{\genus}{genus}
\DeclareMathOperator{\NC}{NC}

\begin{document}

\begin{abstract}
We prove an upper bound for characters of the symmetric groups.
Namely, we show that 
there exists a constant $a>0$ with a property that for every Young diagram
$\lambda$
with $n$ boxes, $r(\lambda)$ rows and $c(\lambda)$ columns
$$ \left| \frac{\Tr \rho^{\lambda}(\pi)}{\Tr \rho^{\lambda}(e)} \right|  \leq 
\left[a
\max\left(\frac{r(\lambda)}{n},\frac{c(\lambda)}{n},\frac{|\pi|}{n}
\right)\right]^{|\pi|}, $$
where $|\pi|$ is the minimal number of factors needed to write $\pi\in S_n$ as a product of 
transpositions. 
We also give uniform estimates for the error term in the Vershik-Kerov's and
Biane's character formulas and give a new formula for free cumulants of the
transition measure.
\end{abstract}

\maketitle

\section{Introduction}

\subsection{Normalized characters}
\label{subsec:normalized-characters}
For a Young diagram $\lambda$ having $n$ boxes and a permutation $\pi\in S_l$ (where
$l\leq n$) we define the \emph{normalized character}
\begin{equation} 
\label{eq:sigma}
\Sigma^{\lambda}(\pi)= (n)_{l}\ \chi^{\lambda}(\pi),
\end{equation}
where $(n)_l=n (n-1) \cdots (n-l+1)$ denotes the falling power and where
$$ \chi^{\lambda}(\pi) = \frac{\Tr \rho^{\lambda}(\pi)}{\Tr \rho^{\lambda}(e)}$$
is the character rescaled in such a way that $\chi^{\lambda}(e)=1$.

% Since for any integers $k_1,\dots,k_m\geq 2$ we identify $(k_1,\dots,k_m)$ with
% a permutation in $S_{k_1+\cdots+k_m}$ with the cycles of length $k_1,\dots,k_m$
% we will also use the notation
% $$ \Sigma^{\lambda}(k_1,\dots,k_m)=\
% (n)_{k_1+\cdots+k_m} \ \chi^{\lambda} (k_1,\dots,k_m). $$

\subsection{Short history of the problem}
Unfortunately, the canonical tool for calculating characters, 
the Murnaghan--Nakayama rule, quickly becomes cumbersome and hence intractable 
for computing characters corresponding to large Young diagrams. Nevertheless 
Roichman \cite{Roichman1996} showed that it is possible to use it to find an upper bound for
characters, namely he proved that there exist constants $0<q<1$ and $b>0$ such
that
\begin{equation}
\label{eq:roichman}
 |\chi^{\lambda}(\pi) | \leq 
\left[\max\left(\frac{r(\lambda)}{n},\frac{c(\lambda)}{n},q\right) \right]^{b\ |\supp \pi|},
\end{equation}
where 
$r(\lambda)$, $c(\lambda)$ denote the numbers of rows and columns of $\lambda$ and
$\supp \pi$ denotes the support of a permutation $\pi$ (the set of its non-fixed points).
Inequality \eqref{eq:roichman} is not satisfactory for many practical purposes (such as 
\cite{MooreRussellSniady}) since it provides rather weak estimates in the case
when the Young diagram $\lambda$ is \emph{balanced}, i.e.\
$r(\lambda),c(\lambda)=O(\sqrt{n})$.
% and $\pi$ is not a very long permutation, i.e.\  $|\supp(\pi)| =
% o(n)$.

Another approach to this problem was initiated by Biane \cite{Biane1998,Biane2003} who 
showed that the value of the normalized character $\Sigma^{\lambda}(\pi)$ 
% (its definition will be recalled in Section \ref{subsec:normalized-characters}) 
can be expressed as a polynomial (called Kerov polynomial) in
\emph{free cumulants} of
the \emph{transition measure} of a Young diagram $\lambda$. The work of Biane
was based on previous
contributions of Kerov \cite{Kerov1993,Kerov1999} and Vershik. Free cumulants of the transition
measure have a nice geometric interpretation therefore Kerov polynomials are a perfect tool 
for study of the character $\chi^{\lambda}(\pi)$ in the limit when the permutation $\pi$ is fixed
and the Young diagram $\lambda$ tends in some sense to infinity.

Unfortunately, despite much progress in this field
(\cite{DolegaF'eray'Sniady2008} and references therein)
our understanding of Kerov polynomials is still not satisfactory; in particular it is not clear
how to use Kerov polynomials in order to obtain non-trivial estimates on the characters
$\chi^{\lambda}(\pi)$ when the length $|\pi|$ of the permutation $\pi\in S_n$ is comparable
with $n$.

In a recent work of one of us with Rattan \cite{Rattan'Sniady2008} we
took yet another approach: thanks to the generalized Frobenius formula we showed
that the value of a normalized character of a given Young diagram $\lambda$ 
can be bounded from above by the value of the normalized character of a rectangular 
Young diagram $p\times q$ for suitably chosen $p,q$. For such a rectangular Young diagram
the value of the normalized character can be explicitly calculated thanks to
the formula of Stanley \cite{Stanley2003/04}. %(which will be recalled as Theorem \ref{theo:stanley} in Section \ref{sec:stanley-feray}).
In this way we proved that for each $C$ there exists a constant $D$ with a property that 
if $r(\lambda),c(\lambda)< C \sqrt{n}$ then
\begin{equation}
\label{eq:rattan-sniady}
|\chi^{\lambda}(\pi)|< \left( \frac{D \max( 1, \frac{|\pi|^2}{n})}{\sqrt{n}} \right)^{|\pi|}, 
\end{equation}
where $|\pi|$ denotes the minimal number of factors necessary to write $\pi$ as a product of 
transpositions. Inequality \eqref{eq:rattan-sniady} gives a much better estimate 
than \eqref{eq:roichman} for balanced Young diagrams and a quite short
permutation ($|\pi|=o(\sqrt{n})$) but a careful analysis of its proof shows
that it gives non-trivial estimates only if 
$\max\big(r(\lambda),c(\lambda)\big) < O(n^{3/4})$.

% For the sake of completeness we point out that the studies of Kerov
% polynomials and of
% the Stanley character formula are very much related to each other
% \cite{Biane2003}.

\subsection{The main result}
% Recently Stanley \cite{Stanley-preprint} conjectured a formula for normalized characters
% of symmetric groups associated with an arbitrary Young diagram 
% which generalized his previous formula for rectangular Young diagrams \cite{Stanley2003/04}.
% This conjecture was proved by F\'eray \cite{F'eray-preprint}. In this article we
% study implications of this formula (which will be referred to as \emph{Stanley-F\'eray formula})
% to the asymptotic estimates of the characters.
Our main result is the following inequality.
\begin{theorem}
\label{theo:rough}
There exists a constant $a>0$ with a property that for every Young diagram $\lambda$
%  $$ |\Sigma^{\lambda}_{\mu}| \leq n^{|\supp \mu|} (a \balanced)^{|\mu|}, $$
 $$ |\chi^{\lambda}(\pi)| \leq 
\left[a
\max\left(\frac{r(\lambda)}{n},\frac{c(\lambda)}{n},\frac{|\pi|}{n}
\right)\right]^{|\pi|}, $$
where $n$ denotes the number of boxes of $\lambda$.
\end{theorem}
It is easy to check that \eqref{eq:rattan-sniady} is a consequence of this theorem
and that it gives better estimates than \eqref{eq:roichman} if $\frac{r(\lambda)}{n},
\frac{c(\lambda)}{n},\frac{|\pi|}{n}$ are smaller than some positive
constant.
It is natural to ask what is the optimal value of the constant $a$. Asymptotics
of characters of symmetric groups related to Thoma characters shows that $a\geq 1$.

% \begin{conjecture} There is a constant $C$ with a property that
% $$ |\chi^{\lambda}(\pi) | \leq 
% C \left[ \max\left(\frac{r(\lambda)}{n},\frac{c(\lambda)}{n}\frac{|\supp(\pi)|}{n}\right)\right]^{|\pi|}.$$
% \end{conjecture}

\subsection{Young diagrams}
In the following we shall identify a Young diagram $\lambda$
with the set of its boxes which we regard as a subset of $\N^2$
given by a graphical representation of $\lambda$ according to the French notation;
namely, for a partition $\lambda=(\lambda_1,\dots,\lambda_k)$ 
it is the set
\begin{equation}
\label{eq:identifyNatural}
\lambda = \bigcup_{1\leq i\leq k}   \{1,2,\dots,\lambda_i\} \times  \{i\}=
\{(p,q)\in\N^2: 1\leq p \leq \lambda_q\}.   
\end{equation}

\subsection{The main tool: reformulation of Stanley character formula}
\label{subsec:mainresult}

Our main tool in our investigations will be the following reformulation of 
Stanley character formula.

The set of cycles of a permutation $\pi$ is denoted by $C(\pi)$.
For given permutations $\sigma_1,\sigma_2\in S_l$ we shall consider 
colorings $h$ of the cycles of $\sigma_1$ (where each cycle is colored by the
number of some column of $\lambda$) and of the cycles of $\sigma_2$ 
(where each cycle is colored by the number of some row of $\lambda$). Formally,
each such coloring can be viewed as a function
$h:C(\sigma_1) \sqcup C(\sigma_2)\rightarrow\N$.
%  as a pair of functions
% $h_1,h_2:\{1,\dots,l\}\rightarrow \N$ with a property that $h_1$ is $\sigma_1$-invariant
% and $h_2$ is $\sigma_2$-invariant or 
% as a 
% For given permutations $\sigma_1,\sigma_2\in S_l$ 
We say that a coloring $h$
%  $h:C(\sigma_1) \sqcup C(\sigma_2)\rightarrow\N$ 
is compatible with a Young diagram $\lambda$ if for all
$c_1\in C(\sigma_1)$ and $c_2\in C(\sigma_2)$ if $c_1\cap c_2\neq \emptyset $ then
$(h(c_1),h(c_2))\in\lambda$; in other words
\begin{equation}
\label{eq:warunek}
0<h(c_1)\leq \lambda_{h(c_2)}
\end{equation}
holds true for all $c_1\in C(\sigma_1)$, $c_2\in C(\sigma_2)$ such that $c_1\cap
c_2\neq\emptyset$.

% The following theorem will be our main tool in the analysis of asymptotics of representations
% of symmetric groups.
\begin{theorem}[The new formulation of Stanley character formula]
\label{theo:character-sniady}
For any Young diagram $\lambda$ and a permutation $\pi\in S_l$ (where $l\leq n$)
the value of the normalized character \eqref{eq:sigma} is given by
\begin{equation}
\label{eq:main-theorem}
 \Sigma^{\lambda}(\pi)= \sum_{\substack{\sigma_1,\sigma_2\in S_l,\\ \sigma_1 \sigma_2=\pi}}
(-1)^{|\sigma_1|}\ N^{\lambda}(\sigma_1,\sigma_2),
\end{equation}
where 
% ${\lambda}(\sigma_1,\sigma_2)$ denotes the number of functions compatible with
% $\sigma_1,\sigma_2$. 
\begin{multline}
\label{eq:usual}
N^{\lambda}(\sigma_1,\sigma_2)= 
\#\{h: h \text{ is a coloring of the cycles of $\sigma_1$ and $\sigma_2$} \\
\text{which is compatible with $\lambda$}\}. 
\end{multline}
\end{theorem}

\begin{example}
\label{example1}
For a given factorization $\pi=\sigma_1 \sigma_2$ it is convenient to consider a 
bipartite graph with the set of vertices $C(\sigma_1)\sqcup C(\sigma_2)$ and with an
edge between vertices $c_1\in C(\sigma_1)$ and $c_2\in C(\sigma_2)$ if and only if
$c_1\cap c_2\neq \emptyset$. Notice that the value of $N^{\lambda}(\sigma_1,\sigma_2)$
does not depend on the exact form of $\sigma_1$ and $\sigma_2$ but only on the corresponding
bipartite graph.

Figure \ref{fig:bipartite} presents such a bipartite graph for
$\pi=(12)$, $\sigma_1=(1)(2)$ and $\sigma_2=(12)$.
Now it becomes clear that 
\begin{align*}
 N^{\lambda}\big( (1)(2), (12) \big) & = \sum_i (\lambda_i)^2.\\
\intertext{Similarly, }
 N^{\lambda}\big( (12), (1)(2) \big) & = \sum_i (\lambda'_i)^2,
\end{align*}
where $\lambda'$ denotes the Young diagram conjugate to $\lambda$. In this way,
Theorem \ref{theo:character-sniady} shows that
\begin{multline*} \Sigma^{\lambda}(12) = n(n-1) \frac{\Tr \rho^{\lambda}(12)}{\Tr \rho^{\lambda}(e)}=\\
N^{\lambda}\big( (1)(2), (12) \big)-N^{\lambda}\big( (12), (1)(2) \big)=
\sum_i (\lambda_i)^2-\sum_i (\lambda'_i)^2.
\end{multline*}

 \begin{figure}
%  \includegraphics{figure02.eps}
%  \caption{Bipartite graph associated to the factorization $(123)=  (13)(2) \cdot (12)(3)$.}
 \includegraphics{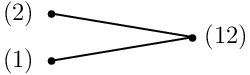}
 \caption{Bipartite graph associated to the factorization $(12)=  (1)(2) \cdot (12)$.}
\label{fig:bipartite}
\end{figure}

%  \begin{figure}
%  \includegraphics{figure00.eps}
%  \caption{ssss}
% \end{figure}

\end{example}

\subsection{Overview of the paper}

% \textbf{It is easy now to compute free cumulants}

In Section \ref{sec:stanley-feray} we will prove the new formulation of
Stanley-F\'eray character formula, Theorem
\ref{theo:character-sniady}.

In Section \ref{sec:random-matrices} we present a relation between the 
characters of symmetric groups and characters of some Gaussian random matrices.
We also give a new formula for calculating free cumulants of (the transition
measure of) a Young diagram. 

Section \ref{sec:technical} is devoted to the proofs of some techical
inequalities.

In Section \ref{sec:asymptotics-characters}
we prove estimates for the characters of the symmetric groups based on
Stanley character formula.

% In Section \textbf{fixme}
% % \ref{sec:leading-term} 
% we study asymptotics of characters of
% symmetric groups for balanced Young diagram and for diagrams with big rows
% and/or columns.

% % We postpone its proof to Section \ref{sec:proof-stanley}.
% The above theorem is a reformulation of Stanley-F\'eray character formula 
% (which will be recalled below in Section \ref{subsec:stanley-feray}), however
% this formulation is more suitable for asymptotic questions 
% (which will be studied in Section \ref{sec:asymptotic}) and it allows 
% a new simple proof of Stanley-F\'eray formula (in Section \ref{sec:new-proof}). 
% 

 \section{Stanley character formula}
 \label{sec:stanley-feray}

In this section, we prove Theorem \ref{theo:character-sniady}. It is quite easy
to show that it is equivalent to a recent formula, conjectured by Stanley
\cite{Stanley-preprint} and proved by the first author \cite{F'eray-preprint}.
But the formulation given here is more useful for the purposes of character
estimates and its proof is more elementary than the one given in
\cite{F'eray-preprint}.

\subsection{Young symmetrizer}
Let $\lambda$ be a Young diagram consisting of $n$ boxes.
In the following we will distinguish the symmetric group $S_n$ which permutes the 
elements $\{1,\dots,n\}$ and the symmetric group $\Sy_n$ which permutes the boxes of $\lambda$.

% We will identify a box of a Young diagram with a pair $(x,y)$, where $x\in\N$  
% is the column and $y\in\N$ is the row of a given box. 
For a box $\Box\in\lambda$ we denote by $r(\Box)\in\N$ (respectively, $c(\Box)\in\N$) the row
(respectively, the column) of $\Box$; in this way
$\Box=(c(\Box),r(\Box))$.

If $\sigma\in \tilde{S}_n$ 
has a property that if two different boxes $\Box_1,\Box_2$ 
are in the same row
then their images $\sigma(\Box_1),\sigma(\Box_2)$ are not in the same column 
% and
% if boxes $\Box_1,\Box_2$ are in the same column
% then $\sigma(\Box_1),\sigma(\Box_2)$ are not in the same row 
then
we define its number of column inversions $\inv(\sigma)$
as the number of pairs $\Box_1,\Box_2$ such that
$\sigma(\Box_1),\sigma(\Box_2)$ are in the same column, $r(\Box_1)<r(\Box_2)$ and
$r(\sigma(\Box_1))>r(\sigma(\Box_2))$. If $\sigma\in \tilde{S}_n$ does not have
this property, then we define $(-1)^{\inv (\sigma)}=0$.

The following theorem gives a very esthetically appealing formula for the
characters of the symmetric groups.

\begin{theorem}
\label{theo:bubu}
Let a Young diagram $\lambda$ having $n$ boxes and $\pi\in S_n$ be given. Let $\hat{\pi}\in \tilde{S}_n$ be a random permutation
distributed with the
uniform distribution on the conjugacy class defined by $\pi$. Then
$$ \chi^{\lambda}(\pi)= \E [(-1)^{\inv(\hat{\pi})}].$$
\end{theorem}

\begin{proof}
We denote 
\begin{align*} 
P_{\lambda}=& \{\sigma\in \Sy_n \colon \sigma \text{ preserves each row
of }\lambda \}, \\
Q_{\lambda}=& \{\sigma\in \Sy_n \colon \sigma \text{ preserves each column
of }\lambda \}
\end{align*}
and define
\begin{align*} 
a_{\lambda} = & \sum_{\sigma\in P_{\lambda}} \sigma \in \C[\Sy_n],\\
b_{\lambda} = & \sum_{\sigma\in Q_{\lambda}} (-1)^{|\sigma|} \sigma \in \C[\Sy_n],\\
c_{\lambda} = &  b_{\lambda} a_{\lambda}.
\end{align*}

It is well-known that $p_{\lambda}=\alpha_{\lambda} c_{\lambda}$ is an idempotent
for a constant $\alpha_\lambda$ which will be specified later. Its image 
$V_{\lambda} = \C[\tilde{S_n}] p_{\lambda}$ under multiplication from the right
on the regular representation  gives a representation
$\rho^{\lambda}$ (where the symmetric group acts by left multiplication)
 associated to a Young diagram $\lambda$. It turns out that
$\alpha_{\lambda}=\frac{\dimm V_{\lambda}}{n!}$.
% We denote by $\rho^L$ the left-regular representation of $\Sy_n$; 
It follows that
for $\tilde{\pi}\in \Sy_n$
\begin{multline}
\label{eq:charakter-babab} 
n! \frac{\Tr \rho^{\lambda}(\tilde{\pi})}{\dimm V_{\lambda}}
= \frac{\Tr \rho^{\lambda}(\tilde{\pi}^{-1})}{\alpha_{\lambda}} = 
%  \Tr \rho^{L}\left(  \tilde{\pi}^{-1}  c_{\lambda}\right) = 
\frac{1}{\alpha_{\lambda}} \sum_{\mu\in \Sy_n} 
 \langle \delta_{\mu}, \tilde{\pi}^{-1} \delta_\mu p_{\lambda}\rangle
= \\
\sum_{\mu\in \Sy_n} \sum_{\tilde{\sigma}_1\in Q_{\lambda}} 
\sum_{\tilde{\sigma}_2\in P_{\lambda}}
(-1)^{|\tilde{\sigma}_1|} [\mu =\tilde{\pi}^{-1} \mu \tilde{\sigma}_1 \tilde{\sigma}_2].
% \sum_{\tilde{\sigma}_1\in Q_{\lambda}} \sum_{\tilde{\sigma}_2\in P_{\lambda}}
% (-1)^{|\sigma_1|} \sum_{\mu\in \Sy_n} 
% [\tilde{\pi}= (\mu^{-1} \tilde{\sigma}_1\mu) (\mu^{-1} \tilde{\sigma}_2 \mu)].
\end{multline}

We define $\hat{\pi}=\mu^{-1} \tilde{\pi} \mu$.
For such a permutation $\hat{\pi}$ there exists at most one factorization
$\hat{\pi}=\tilde{\sigma}_1 \tilde{\sigma}_2$; by Young lemma this factorization exists
if and only if $(-1)^{\inv(\hat{\pi})}\neq 0$. Furthermore, if such a factorization
exists then $(-1)^{|\sigma_1|}=(-1)^{\inv(\hat{\pi})}$. It follows that
$$ n! \frac{\Tr \rho^{\lambda}(\tilde{\pi})}{\dimm V_{\lambda}}=
\sum_{\mu\in \Sy_n} 
(-1)^{\inv(\mu^{-1} \tilde{\pi} \mu)}. $$
As $\mu$ runs over all permutations,
$\hat{\pi}=\mu^{-1} \tilde{\pi} \mu$ runs over all elements of the conjugacy class defined by $\pi$
which finishes the proof.
% 
% Denote $f=\mu^{-1}\circ g$; as $\mu$ runs over all permutations in $\Sy_n$ the corresponding
% $f$ runs over all bijections from $\{1,\dots,n\}$ to the set of boxes of $\lambda$. 
% 
\end{proof}

For permutations $\sigma_1,\sigma_2\in S_l$ we define
$\tilde{N}^{\lambda}(\sigma_1,\sigma_2)$ as the number of one-to-one functions $f$ from 
$\{1,\dots,l\}$ to the set of boxes of $\lambda$ with a property that $r\circ f$ is constant 
on each cycle of $\sigma_2$ and $c\circ f$ is constant on each cycle of $\sigma_1$.

\begin{proposition}
\label{prop:bubu2}
Let $\lambda$ be a Young diagram having $n$ boxes.
For any permutation $\pi\in S_l$ (where $l\leq n$)
$$ \Sigma^{\lambda}(\pi)= \sum_{\substack{\sigma_1,\sigma_2\in S_l\\ \sigma_1 \sigma_2=\pi}}
(-1)^{|\sigma_1|} \tilde{N}^{\lambda}(\sigma_1,\sigma_2).$$
\end{proposition}
\begin{proof}
Let us consider the case $l=n$. Let $\tilde{\pi}\in\tilde{S}_n$ be any 
permutation with the same cycle structure as $\pi\in S_n$.
Our starting point is the analysis of \eqref{eq:charakter-babab}.
Notice that the multliset of the values of $\mu^{-1} \tilde{\pi} \mu$ (over
$\mu\in \Sy_n$)
coincides with the multiset of the values of $f \circ \pi \circ f^{-1}$ (over bijections $f$). 
We define $\sigma_i=f^{-1} \circ \tilde{\sigma_i} \circ f$; then condition
$\mu =\tilde{\pi}^{-1} \mu \tilde{\sigma}_1 \tilde{\sigma}_2$ is equivalent to 
$\pi=\sigma_1 \sigma_2$. It is easy to check that $\tilde{\sigma}_1\in Q_{\lambda}$ if
and only if $c\circ f$ is constant on each cycle of $\sigma_1$ and 
$\tilde{\sigma}_2\in P_{\lambda}$ if
and only if $r\circ f$ is constant on each cycle of $\sigma_2$. Thus
\begin{equation}
\label{eq:ppp}
 n! \frac{\Tr \rho^{\lambda}(\tilde{\pi})}{\dimm V_{\lambda}}
= \sum_{\substack{\sigma_1,\sigma_2\in S_n,\\ \sigma_1\sigma_2=\pi}} 
(-1)^{|\sigma_1|} \tilde{N}^{\lambda}(\sigma_1,\sigma_2) 
\end{equation}
and the proof in the case when $l=n$ is finished.

For a permutation $\sigma\in S_n$ we denote by $\supp \sigma\subseteq\{1,\dots,n\}$ the 
support of a permutation (the set of non-fixed points). We claim that a factorization 
$\pi=\sigma_1 \sigma_2$ has a non-zero contribution to \eqref{eq:ppp} only if
$\supp \sigma_1,\supp \sigma_2\subseteq \supp \pi$. Indeed, if $m\in \supp \sigma_1\setminus
\supp \pi=\supp \sigma_2\setminus \supp \pi$ then for any bijection $f$ at least one of the 
following conditions hold true: $r(f(m))\neq r(f(\sigma_2(m)))$ (in this case
$r\circ f$ is not constant on the cycles of $\sigma_2$) or
$c(f(m))\neq c(f(\sigma_2(m)))$ (in this case $c(f(\sigma_1(\sigma_2(m))))\neq c(f(\sigma_2(m)))$
hence $c\circ f$ is not constant on the cycles of $\sigma_1$).

It follows that if $\pi\in S_l$ then we may restrict the sum in \eqref{eq:ppp} to
factorizations $\pi=\sigma_1 \sigma_2$, where $\sigma_1,\sigma_2\in S_l$.
It remains to notice that 
$$ \tilde{N}_{S_n}^{\lambda}(\sigma_1,\sigma_2) =
(n-l)! \ \tilde{N}_{S_l}^{\lambda}(\sigma_1,\sigma_2), $$
where the quantity on the left-hand side regards $\sigma_1,\sigma_2$ as
elements of $S_n$ while on the quantity on the right-hand side regards them as
elements of $S_l$ and that the analogous relation holds between
$\Sigma^\lambda(\pi)$ for $\pi\in S_n$ and $\pi\in S_l$.
\end{proof}

% \item For any permutation $\pi\in S_l$
% $$ \Sigma^{\lambda}(\pi)= \sum_{\substack{\sigma_1,\sigma_2\in S_l\\ \pi=\sigma_1 \sigma_2 }} 
% (-1)^{|\sigma_2|} \sum_{f} 1, $$
% where the sum runs over all admissible one-to-one functions $f$ from $\{1,\dots,l\}$ to the set of boxes of $\lambda$ . 

% \begin{proposition}
% For any permutation $\pi\in S_l$ (where $l\leq n$)
% $$ \Sigma^{\lambda}(\pi)= \sum_{\substack{\sigma_1,\sigma_2\in S_l\\ \sigma_1 \sigma_2=\pi}}
% {\lambda}(\sigma_1,\sigma_2),$$
% where ${\lambda}$ denotes the number of functions $f$ from 
% $\{1,\dots,l\}$ to the set of boxes of $\lambda$ with a property that $r\circ f$ is constant 
% on each cycle of $\sigma_2$ and $c\circ f$ is constant on each cycle of $\sigma_1$.
% \end{proposition}
\subsection{Forgetting injectivity}

For a pair of permutations $\sigma_1,\sigma_2\in S_l$ we define 
$\hat{N}^{\lambda}(\sigma_1,\sigma_2)$ 
% (respectively, $\hat{N}^{\lambda}(\sigma_1,\sigma_2)$)
as the number of all functions 
% (respectively, all one-to-one functions)
$f:\{1,\dots,l\}\rightarrow \lambda$ (with values in the set of boxes of $\lambda$)
with a property that $r\circ f$ is constant 
on each cycle of $\sigma_2$ and $c\circ f$ is constant on each cycle of $\sigma_1$.

\begin{lemma}
\label{lem:kasowanie}
For any Young diagram $\lambda$ and a permutation $\pi\in S_l$
\begin{equation}
\label{eq:stanley-feray-sniady-2}
\sum_{\substack{\sigma_1,\sigma_2\in S_l\\ \sigma_1 \sigma_2=\pi,}}
  (-1)^{|\sigma_1|} \tilde{N}^{\lambda}(\sigma_1,\sigma_2)=
\sum_{\substack{\sigma_1,\sigma_2\in S_l\\ \sigma_1 \sigma_2=\pi,}}
  (-1)^{|\sigma_1|} \hat{N}^{\lambda}(\sigma_1,\sigma_2).
\end{equation}
\end{lemma}
\begin{proof}
For a given function $f:\{1,\dots,l\}\rightarrow \lambda$ we will show that it
has the same contribution to the left-hand side and to the right-hand side.
Clearly, it is enough to consider the case when $f$ is not a one-to-one
function. It follows that there exists a transposition $\mu\in S_l$ with a
property that $f$ is
constant on the orbits of $\mu$.
Function $f$ contributes to the right-hand side
%$$\sum_{\substack{\sigma_1,\sigma_2\in S_l\\ \sigma_1 \sigma_2=\pi,}}
%  (-1)^{|\sigma_1|} \tilde{N}^{\lambda}(\sigma_1,\sigma_2) $$
with multiplicity
\begin{equation}
\label{eq:contribution}
 \sum_{\substack{\sigma_1,\sigma_2\in S_l\\ \sigma_1 \sigma_2=\pi}} (-1)^{|\sigma_1|},
\end{equation}
where the sum runs over pairs $(\sigma_1,\sigma_2)$ with a property that $\sigma_1 \sigma_2=\pi$
and
$c\circ f$ is constant on each cycle of $\sigma_1$ and 
$r\circ f$ is constant on each cycle of $\sigma_2$.
Map $(\sigma_1,\sigma_2)\mapsto (\sigma_1', \sigma_2')$ with $\sigma'_1=\sigma_1
\mu$,
$\sigma_2'=\mu \sigma_2$ is an involution of
the pairs $(\sigma_1,\sigma_2)$ which contribute to \eqref{eq:contribution};
the only less trivial condition which should be verified is that
$c\circ f$ is constant on each cycle of $\sigma'_1$ but this is equivalent to
$c\circ f$ being constant on each cycle of $\sigma'_1{}^{-1}=\mu \sigma_1^{-1}$.

Since $(-1)^{|\sigma_1|}= (-1)\cdot (-1)^{|\sigma'_1|}$ therefore the contributions
of the pairs $(\sigma_1,\sigma_2)$ and $(\sigma_1',\sigma_2')$
to \eqref{eq:contribution} cancel. In this way we proved that \eqref{eq:contribution} is 
equal to zero which finishes the proof.
% $$ \sum_{\substack{\sigma_1,\sigma_2\in S_l\\ \sigma_1 \sigma_2=\pi,}}
%   (-1)^{|\sigma_1|} \tilde{N}^{\lambda}(\sigma_1,\sigma_2)=
% \sum_{\substack{\sigma_1,\sigma_2\in S_l\\ \sigma_1 \sigma_2=\pi,}}
%   (-1)^{|\sigma_1|} \hat{N}^{\lambda}(\sigma_1,\sigma_2).$$
\end{proof}

\begin{proof}[Proof of Theorem \ref{theo:character-sniady}]
Proposition \ref{prop:bubu2} and Lemma \ref{lem:kasowanie} show that 
$$ \Sigma^{\lambda}(\pi)= \sum_{\substack{\sigma_1,\sigma_2\in S_l\\ \sigma_1 \sigma_2=\pi}}
(-1)^{|\sigma_1|} \hat{N}^{\lambda}(\sigma_1,\sigma_2).$$
Now it is enough to notice that $N^{\lambda}(\sigma_1,\sigma_2)=\hat{N}^{\lambda}(\sigma_1,\sigma_2)$;
the desired bijection is defined as follows: if $m\in\{1,\dots,l\}$ fulfills
$m\in c_1\cap c_2$ for $c_i\in C(\sigma_i)$ we set $f(m)=(h(c_1),h(c_2))$.
\end{proof}

\subsection{Generalization to Young diagrams on $\R_+^2$}
We may identify a Young diagram  
with a subset of $\R^2$ given by a graphical representation of $\lambda$ 
(according to the French notation). For example, for a partition
 $\lambda=(\lambda_1,\dots,\lambda_k)$  it is the set 
\begin{equation}
\label{eq:identify}
\bigcup_{1\leq i\leq k}   [0,\lambda_i] \times  [i-1,i] .   
\end{equation}
In this way we may consider colorings $h$ of the cycles of permutations
$\sigma_1$ and $\sigma_2$ which take real values instead of natural numbers. 
If we fix some numbering of the cycles in $C=C(\sigma_1)\sqcup C(\sigma_2)$ 
then any such coloring $h\colon C\rightarrow\R_+$ 
can be identified with an element of $\R_+^{|C|}$.

We define
\begin{equation}
\label{eq:generalized}
N^{\lambda}(\sigma_1,\sigma_2)= \vol \{h\in \R_+^{|C|}:h \text{ compatible with } \lambda \}. 
\end{equation}
Notice that in the case when $\lambda\subset\R^2$ is as prescribed by \eqref{eq:identify},
the set of functions $h\in \R_+^{|C|}$ compatible with $\lambda$ is a polyhedron hence there
is no difficulty in defining its volume; furthermore definitions \eqref{eq:usual} and
\eqref{eq:generalized} give the same value.

% In the case of the usual Young diagrams 
% this definition of $N^{\lambda}(\sigma_1,\sigma_2)$ coincides with the previous one 
% \eqref{eq:definicja-lambdy}.
The advantage of the definition \eqref{eq:generalized} is that it allows to define characters
for any bounded set $\lambda\in\R_+^2$, in particular for \emph{generalized Young diagrams} (see \cite{Kerov1993}).
%multi-rectangular Young diagrams $\mathbf{p}\times \mathbf{q}$ for general sequences 
%$p_1,\dots,p_k\geq 0$, $q_1,\dots,q_k\geq 0$ which do not have to be natural numbers---just 
%like in the original papers \cite{Stanley2003/04,Stanley-preprint,F'eray-preprint}.

% \begin{figure}
% \includegraphics{skew.eps} 
% \caption{Example of a skew Young diagram for which Theorem
% \ref{theo:character-sniady}
% is false.}
% \label{fig:skew}
% \end{figure} 

It is very natural therefore to ask if Theorem \ref{theo:character-sniady} holds true
also for skew Young diagrams. Unfortunately, this is not the case as it can be seen
for the skew Young diagram $\lambda=(3,2) \setminus (1)$.

\section{Characters of symmetric groups, random matrices and free probability}
\label{sec:random-matrices}

\subsection{Stanley character formula and random matrices}
Let $\lambda$ be a  Young diagram and $N\geq r(\lambda),c(\lambda)$. We consider a random matrix
$T_\lambda=(t_{ij})_{1\leq i,j\leq N}$ such that
\begin{itemize}
\item its entries $(t_{ij})$ are independent random variables;
\item if $(i,j) \in \lambda$ then $t_{ij}$ is a complex centered Gaussian
variable, that is to say that
$$
\E(t_{ij}) =0, \quad 
\E(t_{ij} \overline{t_{ij}}) = 1, \quad
\E(t_{ij}^2) =0; 
$$
\item otherwise, if $(i,j) \notin \lambda$ then $t_{ij}=0$.
\end{itemize}
%We are interested in this matrix because 
The moments of $T_\lambda
T_\lambda^\star$ are given by a formula which is very similar to 
the Stanley formula for characters (Theorem
\ref{theo:character-sniady}):
\begin{theorem}
\label{theo:random-matrices}
With the definitions above and $\pi \in S_l$ with a cycle decomposition 
$k_1,\ldots,k_r$
\begin{equation}
\label{eq:moments-matrix}
\E\big(\Tr(T_\lambda T_\lambda^\star)^{k_1} \cdots \Tr(T_\lambda
T_\lambda^\star)^{k_r} \big)=\sum_{\substack{\sigma_1,\sigma_2 \in S_l \\
\sigma_1 \sigma_2=\pi}} N^\lambda(\sigma_1,\sigma_2).
\end{equation}
\end{theorem}
\begin{proof}
The first step is to expand the product and the trace on the left-hand side:
\begin{eqnarray*}
\E\big(\Tr(T_\lambda T_\lambda^\star)^{k_1} \ldots \Tr(T_\lambda T_\lambda^\star)^{k_r} \big) \hspace{-5cm}\\
&=&\E\Bigg[ \Big(\sum_{i^1_1,j^1_1,\ldots,i^1_{k_1},j^1_{k_1}} t_{i^1_1 j^1_1}
\overline{t_{i^1_2 j^1_1}} \ldots t_{i^1_{k_1} j^1_{k_1}} \overline{t_{i^1_1
j^1_{k_1}}} \Big)\\
& &\hspace{1,5cm} \ldots \Big(\sum_{i^r_1,j^r_1,\ldots,i^r_{k_r},j^r_{k_r}}
t_{i^r_1 j^r_1} \overline{t_{i^r_2 j^r_1}} \ldots t_{i^r_{k_r} j^r_{k_r}}
\overline{t_{i^r_1 j^r_{k_r}}} \Big) \Bigg]\\
&= &
% &\E\left[\sum_{i_1,j_1,\ldots,i_K,j_K} \prod_{l=1}^K t_{i_l j_l}
% \overline{t_{i_{\pi(l)} j_l}} \right]=
\sum_{i_1,j_1,\ldots,i_l,j_l} \E\left[
\prod_{m=1}^l t_{i_m j_m} \overline{t_{i_{\pi(m)} j_m}} \right].
\end{eqnarray*}
Since random variables $(t_{ij})$ are Gaussian, we can apply Wick formula (see
\cite{Zvonkin1997}) to each summand; in order to do this we need to consider
all ways of pairing factors $(t_{i_m,j_m})$ with factors
$(\overline{t_{i_{\pi(m)} j_m}})$. Each such a pairing can be identified with a
permutation $\sigma\in S_l$: therefore
% $\E(t_{ij}t_{kl})=0$ for all values of $i,j,k,l$, the contribution of a
% pairing can be non-zero only if the pairs are formed by one monomial
% $t_{\ldots}$ and one monomial $\overline{t_{\ldots}}$. As there are $K$
% monomials of each type, such pairings are in bijection with $S(K)$. So, if we
% fixed values of $i_1,j_1,\ldots,i_K,j_K$, we have :
$$\E\left[ \prod_{m=1}^l t_{i_m j_m} \overline{t_{i_{\pi(m)}
j_m}} \right]=\sum_{\sigma \in S_l} \left[ \prod_{m=1}^l
\E\left(t_{i_{\sigma(m)} j_{\sigma(m)}} \overline{t_{i_{\pi(m)} j_m}}\right)
\right]$$
% But, from the definition of $T_\lambda$ it comes :
% \E\left(t_{i_{\sigma(l)} j_{\sigma(l)}} \overline{t_{i_{\mu(l)}
% j_l}}\right)=
Using the definition of $T$, this is equal to
$$
\sum_{\sigma \in S_l} \prod_{m=1}^l  [i_{\sigma(m)}=i_{\pi(m)}] \  [j_{\sigma(m)}=j_m] \
[(i_{\sigma(m)},j_{\sigma(m)}) \in \lambda].
$$

If we plug this in our calculation, the left-hand side of \eqref{eq:moments-matrix} is equal to
\begin{multline*}
%\E\big(\Tr(T_\lambda T_\lambda^\star)^{k_1} \ldots \Tr(T_\lambda
%T_\lambda^\star)^{k_r} \big) =\\
% &=&\sum_{i_1,j_1,\ldots,i_K,j_K} \sum_{\sigma \in S(K)} \left( \prod_{l=1}^K
% \delta_{i_{\sigma(l)},i_{\mu(l)}} \delta_{j_{\sigma(l)} j_l}
% \delta_{(i_{\sigma(l)},j_{\sigma(l)}) \in \lambda} \right).\\
\sum_{\sigma \in S_l} \sum_{i_1,j_1,\ldots,i_l,j_l} \left( \prod_{m=1}^l
[i_{\sigma(m)}=i_{\pi(m)}] \  [j_{\sigma(m)}=j_m] \
[(i_{\sigma(m)},j_{\sigma(m)}) \in \lambda] \right).
\end{multline*}
If we denote $\sigma_1=\pi \sigma_2^{-1} $, $\sigma_2=\sigma$ then the sum over
$\sigma$ can be seen as a sum over all $\sigma_1,\sigma_2\in S_l$ such that
$\sigma_1\sigma_2=\pi$. If a sequence $i_1,\dots,i_l$ (respectively,
sequence $j_1,\dots,j_l$) contributes to the above sum then it must be constant
on each cycle of $\sigma_1$ (respecively, each cycle of $\sigma_2$). It follows
that there is a bijective correspondence between sequences
$i_1,j_1,\dots,i_l,j_l$ which contribute to the above sum and colorings of the
cycles of $\sigma_1$ and $\sigma_2$ which are compatible with $\lambda$.
% A family of integers $(i_1,j_1,\ldots,i_K,j_K)$ is exactly the same thing as a
% function from $\{1,\ldots,K\}$ to $\N^2$. The product of the deltas in the sum
% is equal to $1$ if and only if :
% \begin{itemize}
% \item the values of this function are in the Young diagram $\lambda$;
% \item the function $l \mapsto i_l$ (the line of the box associated to $l$) is
% invariant by $\sigma^{-1} \mu$ (for each $l$, $i_{\sigma(l)}=i_{\mu(l)}$);
% \item the function $l \mapsto j_l$ (the column of the box associated to $l$)
% is invariant by $\sigma$ (for each $l$, $j_{\sigma(l)}=j_l$).
% \end{itemize}
% So, if we put $\sigma_1=\sigma^{-1} \mu$ and $\sigma_2=\sigma$, we find
% exactly the equality of the lemma.
\end{proof}

By comparing the above result with Theorem \ref{theo:character-sniady}, we
obtain the following corollary.
\begin{corollary}
Let $\lambda$ be a Young diagram. 
% and $T$ be the corresponding random matrix.
Then for any permutation $\pi\in S_l$ with a cycle decomposition $k_1,\dots,k_r$
 $$|\Sigma^{\lambda}(\pi)| \leq  \E\big(\Tr(T_\lambda T_\lambda^\star)^{k_1}
\cdots \Tr(T_\lambda
T_\lambda^\star)^{k_r} \big).$$
\end{corollary}
We shall not follow this idea in this article and we will prove all estimates
from scratch, but it is worth noticing that the above Corollary shows that
the asymptotics of characters of symmetric groups can be deduced from the
corresponding asymptotics of random matrices. In particular it follows  that 
% Inhis article we are mostly concerned with the situation 
when the lengths $k_1,\dots,k_r$ of the cycles of $\pi$ are big enough then the
asymptotics of the corresponding character is related to the limit distribution
of the largest eigenvalue of $T_{\lambda}T\gwia_{\lambda}$
\cite{SinaiSoshnikov1998}. Notice, however, that due to the minus sign in
Theorem \ref{theo:character-sniady} and the resulting cancelations the 
character $|\Sigma^{\lambda}(\pi)|$ could, at least in priciple,  be much
smaller than the appropriate moment of the random matrix
$T_{\lambda}T_{\lambda}\gwia$.

\subsection{Free cumulants of the transition measure}
\label{subsec:cumulants}

% 
% 
% 
% Since the other parts of this article do not use the results presented in this
% section (except for Theorem
%   \ref{theo:cumulants}) will be quite 
% short in recalling the necessary notions.

For a (continuous) Young diagram $\lambda$ we denote by $\mu^{\lambda}$ its transition measure 
(which is a probability measure on the real line) \cite{Kerov1993,Biane1998} 
and by $R_{m}^{\lambda}:=R_m(\mu^{\lambda})$ we denote the $m$-th free cumulant of $\mu^{\lambda}$.
The importance of free cumulants $R_m^{\lambda}$ in the study of the asymptotics of symmetric groups 
was pointed out by Biane \cite{Biane1998}. 
The following theorem gives a new formula for the free cumulants
$R_m^{\lambda}$. It has a big advantage that it does not involve the notion of
the transition measure and it is related directly with the shape of a Young
diagram.

\begin{theorem}\label{theo:cumulants}
For any Young diagram $\lambda$
\begin{equation} 
\label{eq:freecumulants-come}
R_{l+1}^{\lambda}= \sum_{\substack{\sigma_1,\sigma_2\in S_l,\\
\sigma_1 \sigma_2=(1,2,\dots,l), \\
|\sigma_1|+|\sigma_2|=|(1,2,\dots,l)|}} (-1)^{|\sigma_1|}\
N^{\lambda}(\sigma_1,\sigma_2),
\end{equation}
where the sum runs over minimal factorizations of a cycle of length $l$.
\end{theorem}
\begin{proof}
For a Young diagram $\lambda$ and $c>0$ we denote by $c\lambda$ the
(generalized) Young diagram obtained from $\lambda$ by similarity with scale
$c$. A function $f$ on the set of (generalized) Young diagrams
is said to be homogeneous of degree $m$ if
$$f(c \lambda)= c^m f(\lambda)$$
holds true for all choices of $\lambda$ and $c$. Each free cumulant $R_m$ is homogenous of degree $m$.

%We regard $R_m$ as a function on the set of Young diagrams. 
%The algebra generated by polynomials $(R_{m}: m\in\{2,3,\dots\} )$ can be equipped in a gradation 
%given by $\deg R_{m}^{\lambda}= m$. This gradation is motivated by the equality
%$$ R_m^{c \lambda}= c^m R_m^{\lambda},$$
%where $c\lambda$ denotes the Young diagram obtained from $\lambda$ by similarity with scale $c>0$;
%in other words $R_m$ lambda is homogenous (of degree $m$) with respect to rescaling of the Young diagram.

The value of the normalized character $\Sigma^{\lambda}(1,2,\dots,l)$ on a cycle
can be expressed as a polynomial (known as Kerov polynomial) in free cumulants
$(R^{\lambda}_{m}: m\in\{2,3,\dots\} )$:
$$ \Sigma^{\lambda}(1,2,\dots,l)= R_{l+1}^{\lambda} + (\text{terms of lower
degree}) $$
therefore $R^\lambda_{l+1}$ is the homogeneous part of
$\Sigma^{\lambda}(1,2,\dots,l)$ with degree $l+1$.
We apply \eqref{eq:main-theorem} for $\pi=(1,2,\dots,l)$;
it is easy to see that each summand on the right-hand side is homogeneous of degree
$|C(\sigma_1)|+|C(\sigma_2)|$ which finishes the proof.
\end{proof}

\subsection{Generalized circular operators}
Let $\D$ be the algebra of continuous  functions on $\R_+$. We equip it with an
expected value
$\phi:\D\rightarrow \C$ given by $\phi(f)=\int_0^{\infty} f(t) dt$.

We consider an operator-valued probability space, which by definition is some
$\ast$-algebra $\A$ 
such that $\D\subseteq \A$ and equipped with a conditional expectation $\E:\A\rightarrow\D$.
For a given (generalized) Young diagram $\lambda$ let $T\in\A$ be a generalized
circular operator \cite{VoiculescuDykemaNica1992,Speicher1998}
with a covariance
\begin{align}
\nonumber \big[ k(T,f T\gwia)\big] (s)&= \int_{t:  (t,s)\in\lambda }f(t)\ dt,
\\
\label{eq:strangesign}
\big[ k(T\gwia,f T)\big] (s)&= -\int_{t:  (s,t)\in\lambda }f(t)\ dt, \\
\nonumber \big[ k(T,f T)\big] (s)&= 0,\\
\nonumber \big[ k(T\gwia,f T\gwia)\big] (s)&= 0.
\end{align}

\begin{theorem}
\label{theo:free-proba-may-help}
For any (generalized) Young diagram $\lambda$
$$ R_{l+1}^{\lambda}= \phi\big[ (T\gwia T)^l \big]. $$ 
\end{theorem}
\begin{proof}
It is easy to check that for permutations $\sigma_1,\sigma_2$ which contribute
to \eqref{eq:freecumulants-come} the corresponding bipartite graph is a tree
therefore the calculation of $N^{\lambda}(\sigma_1,\sigma_2)$ is particularly
simple, namely it is a certain iterated integral. The same iterated integral
appears in the nested evaluation of amalgamated free cumulants therefore
$N^{\lambda}(\sigma_1,\sigma_2)=\pm
\phi\big[ k_{\sigma_2}(T\gwia,T,\cdots,T\gwia,T)\big]$.
The plus/minus sign is due to the minus sign in the covariance
\eqref{eq:strangesign}. It is easy to check that in fact
$$(-1)^{|\sigma_1|} N^{\lambda}(\sigma_1,\sigma_2)=
\phi\big[ k_{\sigma_2}(T\gwia,T,\cdots,T\gwia,T)\big].$$

The moment-cumulant formula 
$$\E\big[ (T\gwia T)^l \big]= \sum_{\sigma\in\NC_2}
k_{\sigma}(T\gwia,T,\cdots,T\gwia,T) $$
finishes the proof since there is a bijective correspondence between
non-crossing pair partitions and the minimal factorizations.
% sending $\sigma$ to $\sigma_2$.
\end{proof}

% Except for the minus sign on the right-hand side of \eqref{eq:strangesign} the
% covariance of $T$ looks in the same way as the covariance of generalized
% circular operators corresponding to the limits of Gaussian band matrices
% \cite{Shlyakhtenko1996} which should not be surprising in the light of Theorem
% \ref{theo:random-matrices}. We refer also to the articles
% \cite{'Sniady2003a,DykemaHaagerup2004a} where a similar circular operator is
% considered.

Due to the analytic machinery of free probability the calculation of the
moments of $T$ in a closed form is possible in many cases therefore
Theorem \ref{theo:free-proba-may-help} gives a practical method of calculating
the free cumulants of the Young diagrams.

\section{Technical estimates}
\label{sec:technical}

\subsection{Estimates for the number of colorings
$N^\lambda(\sigma_1,\sigma_2)$}
% \subsection{Estimates for $N^\lambda(\sigma_1,\sigma_2)$}
\label{eq:estimate-lambda}

As we already mentioned in Example \ref{example1} to permutations $\sigma_1,\sigma_2$
we can associate a bipartite graph $C(\sigma_1)\sqcup C(\sigma_2)$ with an edge between
$c_1\in C(\sigma_1)$ and $c_2\in C(\sigma_2)$ if $c_1\cap c_2\neq \emptyset$.

For a bipartite graph $G=C_1\sqcup C_2$ (not necessarily arising from the above construction)
and a Young diagram $\lambda$ we define $N^{\lambda}(G)$ as the number of colorings $h$
of the vertices of $C_1\sqcup C_2$ which are compatible with the Young diagram $\lambda$
(the definition of compatibility in this context is a natural extension of the old one, 
i.e.\ we require that if $c_1\in C_1$ and $c_2\in C_2$ are connected by an edge then
$\big( h(c_1),h(c_2) \big)\in \lambda$).

% In Section \ref{subsec:stanley-feray} we considered sets $C_1,C_2$ and a relation $\sim$
% given by some permutations $\sigma_1,\sigma_2$. Now we are going to consider a more general situation
% where $G=(C,\sim)$ is a finite bipartite graph, where $C=C_1\sqcup C_2$ and if a pair 
% $(c_1,c_2)$ is connected by an edge for $c_1\in C_1$ and $c_2\in C_2$ we denote it by
% $c_1\sim c_2$. Our only assumption is that the degree of every vertex in $C$ is at least one.
% Notice that the definition \eqref{eq:definicja-lambdy} can be easily adapted to this more general 
% case:
% % \begin{equation} 
% % \label{eq:definicja-lambdy}
% $$\lambda(G):=
% \#\{f\colon C\rightarrow \N_+ : f \text{ compatible with }\lambda\}.$$
% % \end{equation}
% Our goal in this section is to find reasonable estimates for $\lambda(G)$.

% We begin with analysis of some special graphs.
We denote by $G_{p,q}$ a full bipartite graph for which $|C_1|=p$ and $|C_2|=q$.

\begin{lemma}
\label{lem:grafGprim}
% \label{enum:p4} 
Let $G$ be a finite bipartite graph with a property that the degree of any vertex is non-zero.
It is possible (not necessarily in a unique way) to remove some of the edges of
$G$ in such a way that the resulting graph $\Grozciete$
is a disjoint union of the graphs of the form $G_{1,1}$, $G_{k,1}$, $G_{1,k}$.

Assume that a Young diagram $\lambda$ consists of $n$ boxes. Then, for any
$A \geq r(\lambda),c(\lambda)$
\begin{equation} 
\label{eq:szacowanie-lambdy}
N^{\lambda}(G)\leq A^{\text{(number of vertices of $G$)}} 
\left(\frac{n}{A^2}\right)^{\text{(number of connected components of
$\Grozciete$)}}.
\end{equation}
% has the following properties:
% \begin{itemize}
% \item every vertex of $\Grozciete$ has non-zero degree;
% \item 
% \end{itemize}
% \end{enumerate}
\end{lemma}
\begin{proof} 
% In order to proof point (\ref{enum:p4}): 
If the graph $G$ contains an edge which connects two vertices of degree bigger than one we 
remove it and iterate this procedure; if no such edge exists then the resulting graph $\Grozciete$ 
has the desired property.

Clearly, $N^\lambda(G)\leq N^\lambda(\Grozciete)$ therefore it is enough to find a suitable upper bound
for $N^\lambda(\Grozciete)$. Since both sides of \eqref{eq:szacowanie-lambdy}
are multiplicative with
respect to the disjoint sum of graphs it is enough to prove \eqref{eq:szacowanie-lambdy}
for $\Grozciete\in\{G_{1,1},G_{k,1}, G_{1,k}\}$. Now, the lemma follows from:
$${N^\lambda}(G_{k,1})=\sum_i \lambda_i^k\leq  
\sum_i \lambda_i A^{k-1} = n A^{k-1}$$
and the analogous inequality for ${N^\lambda}(G_{1,k})$. 
\end{proof}

\begin{proposition}
\label{prop:bipartite-graphs}
Suppose that $r(\lambda),c(\lambda) \leq A \leq n$ and $\sigma_1,\sigma_2\in S_l$ and $\pi=\sigma_1 \sigma_2$.
Then
\begin{multline}
\label{eq:swider}
N^\lambda(\sigma_1,\sigma_2)  \leq  A^{|C(\sigma_1)|+|C(\sigma_2)|}\left(\frac{n}{A^2}\right)^{\orbits(\sigma_1,\sigma_2)} 
\\ \leq   A^{l-|C(\pi)|} n^{|C(\pi)|}
\left(\frac{1}{A}\right)^{o(\sigma_1,\sigma_2)},
\end{multline}
where $\orbits(\sigma_1,\sigma_2)$ denotes the number of orbits in the action of 
$\langle \sigma_1,\sigma_2\rangle$ on the set
$\{1,\dots,n\}$ and
$$o(\sigma_1,\sigma_2)=l-|C(\sigma_1)|-|C(\sigma_2)| +\orbits(\sigma_1,\sigma_2).$$
\end{proposition}
\begin{proof}
The first inequality is a simple corollary from Lemma \ref{lem:grafGprim} since
$A^2\geq r(\lambda) c(\lambda) \geq n$ and the number of connected components of $\tilde{G}$ is bounded from below
by $\orbits(\sigma_1,\sigma_2)$.
The second inequality follows by multiplying by 
\begin{displaymath} \left( \frac{n}{A}
\right)^{|C(\pi)|-\orbits(\sigma_1,\sigma_2)} 
\geq 1.\qedhere
\end{displaymath}
\end{proof}

\subsection{Estimates for the number of factorizations}
Now, we have to find a bound of the number of factorizations of $\pi$ with a
given value of the statistic 
$o(\sigma_1,\sigma_2)$.
\begin{lemma}
\label{lem:feray}
Let $\pi,\sigma_1,\sigma_2\in S_l$ be such that $\pi=\sigma_1 \sigma_2$. 
% Suppose that the action of $\left\langle \sigma_1,\sigma_2 \right\rangle$ on $\{1,\ldots,k\}$ has $r-s$ orbits.
There exist permutations $\sigma_1',\sigma_2'\in S_l$ such that
$\pi=\sigma_1' \sigma_2'$, $|\sigma_1'|+|\sigma_2'|=|\sigma_1|+|\sigma_2|$,
$|\sigma_2'|= |\sigma_2'\sigma_2^{-1} |+|\sigma_2|$   and every
cycle of $\sigma'_1$ is contained in some cycle of $\sigma_2'$.
Furthermore, $|\sigma_1'| = o(\sigma_1,\sigma_2)$.
\end{lemma}
\begin{proof}
If every cycle of $\sigma_1$ is contained in some cycle of $\sigma_2$ then $\sigma_1'=\sigma_1$ and
$\sigma_2'=\sigma_2$ have the required property.

Otherwise, there exist $a,b\in\{1,\dots,n\}$
such that $a$ and $b$ belong to the same cycle of $\sigma_1$ but not
to the same
cycle of $\sigma_2$. We define $\sigma_1'=\sigma_1 (a,b)$, $\sigma_2'=(a,b) \sigma_2$. Notice that $|\sigma_1'|=|\sigma_1|-1$,
$|\sigma_2'|=|\sigma_2|+1$, and the orbits under the action of the subgroups
$\langle \sigma_1,\sigma_2 \rangle$ and $\langle \sigma_1',\sigma_2' \rangle$
are the same. 
%Indeed, if we look at partitions with the refinement order, $\Pi$ (resp. $\Pi'$) is just $C(\sigma_1) \vee C(\sigma_2)$
%(resp. $C(\sigma_1') \vee C(\sigma_2')$). But we have  
% $C(\sigma_1')  \leq C(\sigma_1) \leq \Pi$ and $C(\sigma_2') \leq \Pi$ because $C(\sigma_2) \leq \Pi$ 
% and $a$ and $b$ are in the same part of $\Pi$. So $\Pi' \leq \Pi$.

We iterate this procedure
if necessary (it will finish after a finite number of steps because the length of $\sigma_1$ decreases in each step). 
It remains to prove that
%$|\sigma_2'| \geq |\sigma_2'\sigma_2^{-1} |+|\sigma_2|$; in order to do this
$|\sigma_2'| \geq |\sigma_2'\sigma_2^{-1} |+|\sigma_2|$ (the opposite
inequality follows from the triangle inequality): notice that $|\sigma_2'|-|\sigma_2|$ is equal to $k$
(where $k$ is the number of steps after which the procedure has terminated) and 
$\sigma_2'\sigma_2^{-1} $ is a product of $k$ transpositions, hence
$|\sigma_2'\sigma_2^{-1}|\leq k $.

Furthermore, as every cycle of $\sigma_1'$ is contained in some cycle of
$\sigma_2'$, the set of orbits of the group $\langle
\sigma_1',\sigma_2'\rangle$ is equal to $C(\sigma_2')$, so 
% $  |C(\sigma_2')|=
% \orbits(\sigma_1',\sigma_2')$ therefore
\begin{displaymath} 
o(\sigma_1,\sigma_2)= o(\sigma_1',\sigma_2')= l-
|C(\sigma_1')|= |\sigma_1'|. \qedhere
\end{displaymath}
% which shows $\kappa(\sigma_2')\leq \kappa(\pi)$.
\end{proof}

\begin{lemma}
\label{lem:nie-za-duzo}
For any integers $l\geq 1$ and $i\geq 0$ and for any $\pi\in S_l$
\begin{equation}
\label{eq:iloscpermutacji}
 \# \{ \sigma\in S_l : | \sigma | = i \} \leq \frac{l^{2i}}{i!}.
\end{equation}
\end{lemma}
\begin{proof}
Since every permutation in $S_l$ appears exactly once in the product
$$ [1+ (12) ] [1+ (13)+(23) ] \cdots [1+(1l)+\cdots+(l-1,l)], $$
we have
$$ \sum_i x^i\ \#\big\{ \sigma\in S_l : | \sigma | = i \big\} = (1+x) (1+2x)
\cdots \big(1+(l- 1) x\big).$$
Each of the coefficients of $x^k$ on the right-hand side is bounded from above
by the
corresponding coefficient of $e^{x} e^{2x} \cdots
e^{(l-1)x}=e^{\frac{l(l-1)x}{2}}$, finishing the proof.
\end{proof}

\begin{lemma}
\label{lem:numer-of-cutted-graphs}
There exists a constant $C_0$ with a property that for any $k$
the number of minimal factorisations $\sigma_1\sigma_2=\pi$,
$|\sigma_1|+|\sigma_2|=|\pi|$ of a cycle $\pi=(1, \ldots ,k)$ and
such that the associated graph $\tilde{G}$
consists of $s\geq 2$ components is bounded from above by 
$$\frac{(C_0 k)^{2s-2}}{(2s-2)!}.$$
\end{lemma}
\begin{proof}
Since the factorization is minimal therefore the graph $G$ associated to
$\sigma_1,\sigma_2$ is a tree. We will give to $G$ a structure of planted planar
tree: the root
is the cycle of $\sigma_1$ containing $1$ and his left-most edge links it
to the cycle of $\sigma_2$ containing $1$.

%\correc{; furthermore it carries a natural structure of 
%a plane tree. We select some root, for example as the vertex which corresponds
%to the cycle of $\sigma_1$ which passes through the element $1$. We remove from
%the tree $G$ all leaves and denote the resulting tree by $\Gobciete$.}{}
In each connected component of $\Grozciete$ there is at most one vertex of
degree higher
than one and we shall decorate this vertex. If in some connected component of
$\Grozciete$ there are no such vertices we decorate any of them which is not a
leaf in $G$. In this way
the decorated vertices can be identified with connected components of
$\Grozciete$.

%\correc{}{We will give to $G$ a structure of planted planar tree : the root is the decorated vertex in the same component of $\Grozciete$ as the cycle $c_1$ of $\sigma_1$ containing $1$ and his left-most edge links it either to the cycle of $\sigma_2$ containing $1$ if the root is $c_1$ or to $c_1$ else.}

We consider the graph $\Gobciete$ obtained from $G$ by removing
the leaves (except the root) and the
graph $\GobcieteBardzo$ which consists of the decorated vertices of
$\Grozciete$; we connect two vertices $A,B\in\GobcieteBardzo$ by an edge if
vertices $A,B$ are connected in $G$ (or, equivalently, $\Gobciete$) by a
\emph{direct path}, i.e.\ a path which does not pass through any
connected component of $\tilde{G}$ other than the ones specified by $A$ and
$B$.
% vertex of
% $\GobcieteBardzo$. 
It is easy to see that $\GobcieteBardzo$ inherits the
structure of a plane rooted tree from $G$ (we define the root of
$\GobcieteBardzo$ to be the connected component of $\Grozciete$ of the root of
$G$) and it has $s$ vertices. It follows that the number of such trees
$\GobcieteBardzo$ is bounded from above by the Catalan number $\frac{1}{s+1}
\binom{2s}{s}<4^s$.

In order to reconstruct the tree $\Gobciete$ from $\GobcieteBardzo$ we have to
specify for each edge of $\GobcieteBardzo$ if it comes from a single edge of
$\Gobciete$ or from a pair or a triple of consecutive edges of $\Gobciete$; it
follows that we have (at most) $3^{s-1}$ choices. It might happen also for two
adjacent (with respect to the planar structure) edges $e_1,e_2$ of
$\GobcieteBardzo$ that each of these edges $e_i$ corresponds to a pair or a
triple of consecutive edges $f_i=(f_{i1},f_{i2}[,f_{i3}])$ of $\Gobciete$ and
these tuples $f_1$ and $f_2$ have one edge in common. There are
at most $2s-3$ such
pairs of adjacent edges which accounts for at most $2^{2s-3}$ choices. If
the root of $G$ is not a decorated vertex, it might happen that it is a leaf
or that it belongs to the left-most and/or to the right-most edge of the
root of $\Gobciete$: there are $4$ choices for that.

In order to reconstruct tree $G$ from $\Gobciete$ we have to specify if the root
of $G$ is a decorated vertex or not. Furthermore we have to specify places in
which we will add missing $l$ leaves to the tree $\Gobciete$ (note that $l\leq
k+1-s$); it is easy to see that this is equivalent to specifying a partition
$l=a_1+\cdots+a_{2s-1}$, where $a_1,\dots,a_{2s-1}\geq 0$ are integers. It
follows that the number of choices is bounded from above by 
$$ 2 \binom{l+2s-2}{2s-2}\leq 2 \binom{k+s-1}{2s-2}\leq 2
% \frac{(k+\frac{1}{2})^{2s-2}}{(2s-2)!}. $$ 
\frac{\left(k+\frac{1}{2}\right)^{2s-2}}{(2s-2)!}. $$
%%%%%%%%%%%%%% UWAGA! Ta nierownosc jest troszke delikatna.

A minimal factorisation is determined by its bicolored map with
a marked edge, for example the one linking the two cycles containing $1$
\cite{GouldenJackson1992}. With our
construction, the coloring is determined by the root which always belongs to $C(\sigma_1)$ and the marked edge is its
left-most edge.

It follows that the total number of choices is bounded from above by
% $$2\cdot 3^{s-1} \cdot 4^s \frac{(k+\frac{1}{2})^{2s-2}}{(2s-2)!}.$$
\begin{displaymath}
2 \cdot 3^{s-1}\cdot 2^{2s-1} \cdot 4^s
\frac{\left(k+\frac{1}{2}\right)^{2s-2}}{(2s-2)!}.\qedhere 
\end{displaymath}
\end{proof}

\section{Asymptotics of characters}
\label{sec:asymptotics-characters}

\subsection{Upper bound for characters: proof of Theorem \ref{theo:rough}}
\begin{proof}[Proof of Theorem \ref{theo:rough}]
Let $k_1,\dots,k_r\geq 2$ be the lengths of the non-trivial cycles in the cycle
decomposition of $\pi\in S_n$.
It follows that $l:=k_1+\cdots+k_r=|\supp \pi|$ and in the following we will regard $\pi$ as an element of
$S_l$. We denote $A=\max(l,r(\lambda),c(\lambda))$.

We consider a map which to a pair $(\sigma_1,\sigma_2)$ associates any pair 
$(\sigma_1',\sigma_2')$ as prescribed by Lemma \ref{lem:feray}. For any 
fixed $\sigma_2'$ the
permutations $\sigma_2$ such that $|\sigma_2'|= |\sigma_2'\sigma_2^{-1} |+|\sigma_2|$
can be identified with non-crossing partitions of the cycles of $\sigma_2'$ 
(see \cite[Section 1.3]{Biane1997}).  It follows that the number of such permutations
$\sigma_2$ is equal to the product of appropriate Catalan numbers and, hence, this product is bounded from above by $4^{l}$. Therefore Theorem \ref{theo:character-sniady} and Proposition \ref{prop:bipartite-graphs} show that
\begin{multline*}
% \label{eq:w-sumie-pan-niewiele-umie}
|\Sigma^{\lambda}(\pi) | \leq
\sum_{\substack{\sigma_1,\sigma_2\in S_l,\\ \sigma_1 \sigma_2=\pi}} A^{l-r} n^r \left(\frac{1}{A}\right)^{o(\sigma_1,\sigma_2)} \leq
4^l \sum_{\substack{\sigma_1',\sigma_2'\in S_l, \\ \sigma'_1 \sigma'_2=\pi}} A^{l-r} n^r \left(\frac{1}{A}\right)^{|\sigma'_1|}\\
\leq 4^l A^{l-r} n^r \sum_{\sigma'_1\in S_l} \left(\frac{1}{A}\right)^{|\sigma'_1|} \leq
4^l A^{l-r} n^r \sum_{i \geq 0} \frac{l^{2i}}{A^i i!},
\end{multline*}
where the last inequality follows from Lemma \ref{lem:nie-za-duzo}. It follows that
\begin{equation}
|\Sigma^{\lambda}(\pi) | \leq  4^l A^{l-r} n^r e^{\frac{l^2}{A}} \leq (4e)^l A^{l-r} n^r. 
\end{equation}
After dividing by $(n)_l \geq \left(\frac{n}{e}\right)^l$, we obtain $|\chi^{\lambda}(\pi) | \leq (4 e^2)^l (A/n)^{l-r}$,
which finishes the proof because $l=|\supp \pi| \leq 2|\pi|$ and $l-r=|\pi|$.
\end{proof}

\subsection{Error term for balanced Young diagrams}

Biane \cite{Biane1998}
proved that if the permutation $\pi$ is fixed, its non-trivial cycle lengths are
equal to $k_1,\dots,k_r\geq 2$, and the Young diagram $\lambda$
is balanced (i.e.\ $r(\lambda),c(\lambda)=O(\sqrt{n})$) then
$$\chi^\lambda(\pi) = \frac{1}{(n)_{|\supp(\pi)|}}
\prod_{i=1}^r R_{k_i+1}(\lambda) + O\left(n^{-\frac{|\pi|+2}{2}}\right).$$

The following theorem gives a uniform estimate for the error term in
Biane's formula
(note that if $\pi$ is fixed and the Young diagram balanced, 
with $\varepsilon = \frac{C}{ \sqrt{ n } }$ and $A=D\ \sqrt{ n }$ we recover the result of Biane).

\begin{theorem}
\label{theo:equiv-cumulants}
There exists a constant $a$ such that, for any $0 < \varepsilon < 1$, any Young diagram $\lambda$ of size $n$ and any permutation $\pi \in S_n$ such that
$|\supp(\pi)|^2 \leq \varepsilon A$ and $r(\lambda),c(\lambda)\leq A\leq n$ we
have:
% \begin{multline*}
$$\left|\chi^\lambda(\pi) - \frac{1}{(n)_{|\supp(\pi)|}}
\prod_{i=1}^r R_{k_i+1}(\lambda)\right| \leq 
% \frac{1}{(n)_{|\supp(\pi)|}}
\left(\varepsilon^2 +
\frac{A}{n} \varepsilon \right) \left(\frac{a A}{n}\right)^{|\pi|},
$$
% \end{multline*}
where the $k_i$ are the lengths of the non-trivial cycles of $\pi$.
\end{theorem}

\begin{proof}
We can assume that $\pi\in S_l$ has no fixpoints.
Using Theorem \ref{theo:character-sniady} and Theorem \ref{theo:cumulants}
together with the fact that any minimal factorisation of $\pi$ is a product
of minimal factorisations of its cycles, we can write:
$$\Sigma^\lambda(\pi) - \prod_{i=1}^r R_{k_i+1}^\lambda =
\sum_{\substack{\sigma_1,\sigma_2 \in S_l \\ \sigma_1 \sigma_2 = \pi \\
|\sigma_1| + |\sigma_2| > |\pi|}} (-1)^{|\sigma_1|}
N^\lambda(\sigma_1,\sigma_2).$$
% where $l=|\supp \pi|$ is the support of $\pi$.
To such a pair $(\sigma_1,\sigma_2)$ of permutations we can associate
one of the pairs of permutations $(\sigma'_1,\sigma'_2)$ given by Lemma
\ref{lem:feray} with $|\sigma'_1| \geq 1$.

Consider separately the case $|\sigma'_1|=1$. Then $\orbits(\sigma_1,\sigma_2)=
\orbits(\sigma'_1,\sigma'_2) \geq |C(\pi)|-1 $ and the first inequality in
\eqref{eq:swider} shows that
$$N^\lambda(\sigma_1,\sigma_2)  \leq 
A^{|C(\sigma_1)|+|C(\sigma_2)|}\left(\frac{n}{A^2}\right)^{|C(\pi)|-1} $$
% ^{\orbits(\sigma_1,\sigma_2)} $$
therefore the estimate given by
% 
% 
% 
% $$l-|C(\sigma_1)|+|C(\sigma_2)|=|\sigma_1|+|\sigma_2|-l \geq |\pi| +l+2 =
% |C(\pi)|+2$$
% 
% $$o(\sigma_1,\sigma_2)=l-|C(\sigma_1)|-|C(\sigma_2)|
% +\orbits(\sigma_1,\sigma_2).$$
% 
% 
% We necessarely have $|C(\pi)|-\orbits(\sigma_1,\sigma_2)=1$ (it's a
% non-negative 
% odd number less or equal to 1),
Proposition \ref{prop:bipartite-graphs} can be improved to the following
one:
$$N^\lambda(\sigma_1,\sigma_2) \leq A^{l-|C(\pi)|}n^{|C(\pi)|} \frac{1}{n}.$$

Clearly $A\geq l$ therefore by the same
argument as in the proof of Theorem \ref{theo:rough}, we obtain the inequality
$$\left|\Sigma^\lambda(\pi) - \prod_{i=1}^r R_{k_i+1}(\lambda) \right| \leq 4^l
A^{l-r} n^r \left(\frac{l^2}{n} + \sum_{i \geq 2} \frac{l^{2i}}{A^i
i!}\right).$$
The proof is now finished thanks to the remarks of the previous subsection and the
inequality $\exp(z)-1-z \leq z^2$ for $0<z<1$.
% $$\left|\Sigma^\lambda(\pi) - \prod_{i=1}^r R_{k_i+1}(\lambda) \right| 
% \leq 4^l A^{l-r} n^r \left(\frac{l^2}{n} + \sum_{i \geq 2} \frac{l^{2i}}{A^i
% i!}\right)$$
% % \leq 4^l A^{l-r} n^r \sum_{i \geq 1} \frac{l^{2i}}{A^i i!}
% % \leq 4^l A^{l-r} n^r \frac{2 l^2}{A}$$
% where we used that
% % The proof is now over thanks to the remarks of the previous section and the
% % inequality 
% $\exp(z)-1 \leq 2z$ holds true for $0<z<1$.
% 
% $$\frac{\left|\Sigma^\lambda(\pi) - \prod_{i=1}^r R_{k_i+1}(\lambda)
% \right|}{(n)_l} \leq  \left(\frac{(4e)^2 A}{n}\right)^{|\pi|}  
% 2\epsilon  $$
\end{proof}

\subsection{Characters of symmetric groups related to Thoma characters}

Vershik and Kerov \cite{VershikKerov1981} proved that if 
$\pi$ is a fixed permutation with the lengths of non-trivial
cycles $k_1,\dots,k_r$
then for any Young diagram $\lambda$ with $n$ boxes
% \begin{equation}
% \label{eq:thoma} 
% \chi^{\lambda}(\pi)= \prod_{1\leq i\leq r} 
% \left[  \sum_j \alpha_j^{k_i-1} - \sum_j (-\beta_{j})^{k_i-1} \right] 
% + O(n^{-1}), 
$$\chi^{\lambda}(\pi) =
\prod_{i=1}^r \left[ \sum_j
\alpha_j^{k_i} - \sum_j (-\beta_{j})^{k_i} \right] + O\left(\frac{1}{n}\right) 
$$
% \end{equation}
where $\alpha_j=\frac{\lambda_j}{n}$, $\beta_j=\frac{\lambda'_j}{n}$;
we prefer to write this formula in an equivalent form
\begin{equation}
\label{eq:vershik-kerov2}
\chi^{\lambda}(\pi) = \frac{n^l}{(n)_l}
\prod_{i=1}^r \left[ \sum_j
\alpha_j^{k_i} - \sum_j (-\beta_{j})^{k_i} \right] + O\left(\frac{1}{n}\right). 
\end{equation}

In this section we will prove Theorem \ref{theo:estimate-free} which together
with Theorem \ref{theo:equiv-cumulants} give a uniform estimate for
the error term in the formula \eqref{eq:vershik-kerov2}. In particular, for
$A=n$ and $\varepsilon = \frac{D}{n}$ we recover the result of Vershik and Kerov.

% where the error term can be explicitly bounded thanks to Theorem
% \ref{theo:equiv-cumulants} and  Theorem \ref{theo:estimate-free}.
% \end{theorem}

% In the case where one row or one column is large, we have the following
% estimate
% for the free cumulants:
\begin{theorem}
\label{theo:estimate-free}
There exist constants $a,C>0$  with the following property.
Let $k_1,\dots,k_r$ be positive integers; we denote $k_1+\ldots+k_r=l$.
If $\lambda$ is a Young diagram having $n$ boxes with less than $A$ rows and
columns and such that
$\varepsilon=\frac{(k_1^2+\cdots+k_r^2) n}{A^2} < C$ then
\begin{equation}
\label{eq:freecumulants-estimate}
\left| \frac{\prod_{i=1}^r R_{k_i+1}(\lambda)}{n^l} - \prod_{i=1}^r \left[
\sum_j
\alpha_j^{k_i} - \sum_j (-\beta_{j})^{k_i} \right] \right|
\leq  \varepsilon \left(\frac{A}{n}\right)^{l-r} a^r,
\end{equation}
where $\alpha_j=\frac{\lambda_j}{n}$, $\beta_i=\frac{\lambda'_j}{n}$. 
\end{theorem}

\begin{proof}
Firstly, let us consider the case $r=1$. Note that
\begin{align*}
N^\lambda\big(e,(1,\ldots,k)\big) =\sum_j (n \alpha_j)^k,\quad
N^\lambda\big((1,\ldots,k),e\big) =\sum_j (n \beta_j)^k
\end{align*}
therefore Theorem \ref{theo:cumulants} implies that the left-hand side of
\eqref{eq:freecumulants-estimate} is equal to
$$\Bigg|\frac{1}{n^k} \sum_{\substack{\sigma_1,\sigma_2 \in S_k \backslash \{e\}
\\ \sigma_1 \sigma_2 = \pi \\ |\sigma_1| + |\sigma_2| = |\pi|}}
(-1)^{|\sigma_1|} N^\lambda(\sigma_1,\sigma_2)\Bigg|.$$
For a pair of permutations $\sigma_1,\sigma_2$ which contributes to the above
sum  we consider the bipartite graph $G$ and the graph $\tilde{G}$
given by Lemma \ref{lem:grafGprim}. Clearly, in this case graph $\tilde{G}$ 
has more than one component. 
With Lemma \ref{lem:numer-of-cutted-graphs} and Lemma \ref{lem:grafGprim}
\begin{multline}
\label{eq:estimate-a-la-vk}
\left(\frac{n}{A}\right)^{k-1} \left|
\frac{R^\lambda_{k+1}}{n^k} - \sum_j \big( \alpha_j^{k} - 
% \sum_j
(-\beta_{j})^{k} \big) \right|
\leq \\ \sum_{s \geq 2} \frac{(C_0 k)^{2s-2}}{(2s-2)!}
\left(\frac{n}{A^2}\right)^{s-1} \leq
 \frac{2 C_0^2 k^2 n}{A^2}= 2 C_0^2 \varepsilon,
\end{multline}
where the last inequality holds true if $\frac{2 C_0^2 k^2 n}{A^2}=2 C_0^2
\varepsilon$ is smaller than some positive constant
and the proof is finished in
the case $r=1$.

For the general case, we put $\epsilon_i = \frac{k_i^2
n}{A^2}$. We denote 
\begin{align*}
 X_i&=\frac{1}{a} \left(\frac{n}{A}\right)^{k_i-1}
\frac{R^\lambda_{{k_i}+1}}{n^{k_i}},\\
 Y_i&=\frac{1}{a} \left(\frac{n}{A}\right)^{k_i-1} \sum_j \big( \alpha_j^{k_i} -
(-\beta_{j})^{k_i} \big).
\end{align*}
Let us fix $a>2$. Clearly, $|Y_i|< \frac{2}{a}< 1 $ hence
\eqref{eq:estimate-a-la-vk} shows that $|X_i|<1$ if $\varepsilon$ is smaller
than some positive constant. Telescopic summation
\begin{multline*}
X_1 \cdots X_r - Y_1 \cdots Y_r= X_1 \cdots X_{r-1} (X_r-Y_r)+\\
 X_1 \cdots X_{r-2} (X_{r-1}-Y_{r-1}) Y_r + \cdots 
+(X_1-Y_1) Y_2 \cdots Y_r
\end{multline*}
shows that
\begin{multline*} \frac{1}{a^r} \left(\frac{n}{A}\right)^{l-r}
\left| \frac{\prod_{i=1}^r R^\lambda_{k_i+1}}{n^l} - \prod_{i=1}^r \left[ \sum_j
\alpha_j^{k_i} - \sum_j (-\beta_{j})^{k_i} \right] \right| \leq \\
\frac{2 C_0 (\epsilon_1+\cdots+\epsilon_r)}{a}.\qedhere
\end{multline*}

% $$\big|\ldots\big| \leq  \left(\frac{A}{n}\right)^{l-r} \left( \prod_{i=1}^r
% (1+C' \varepsilon_i)] - 1 \right) \leq  \exp(C' \varepsilon -1)
% \left(\frac{A}{n}\right)^{l-r}$$
\end{proof}

\subsection{Concluding remarks}

In the case where $\pi$ is a fixed permutation, we only need Lemma 11 and we
can avoid most of the technicalities. Therefore, our method gives a
unified, simple way to reprove three important results on asymptotics of
character values on fixed permutations as well as new results: the
intermediate case between balanced diagrams ($A=\Theta(\sqrt{n})$) and
diagrams with long rows and/or columns ($A=\Theta(n)$) has, to our
knowledge, not been studied until now. Moreover, it is interesting to note
that our method can be extended to quite long permutations.

\section*{Acknowledgments}

Research of P{\'S} is supported by the MNiSW research grant P03A 013 30, by the EU Research
Training Network ``QP-Applications", contract HPRN-CT-2002-00279 and  by the EC
Marie Curie Host Fellowship for the Transfer of Knowledge ``Harmonic Analysis,
Nonlinear Analysis and Probability", contract MTKD-CT-2004-013389.

VF has been supported by the Project Polonium, one of Egide's PHC to visit
P\'S. He also would like to thank his advisor Philippe Biane.

\bibliographystyle{alpha}
\bibliography{biblio2007}

\end{document}